\documentclass[12pt]{article}
\usepackage{amssymb}
\newcommand{\mysection}[1]{\section{#1}\setcounter{equation}{0}}

\title{\bf 
Exact finite approximations of average-cost countable\\
Markov Decision Processes}

\author{{\bf Arie Leizarowitz}
\\{\small Faculty of Mathematics}\\{\small Technion, Haifa 32000}\\
{\small Israel} \and {\bf Adam Shwartz}
\\{\small Faculty of Electrical Engineering}\\{\small Technion, Haifa 32000}\\
{\small Israel}}

\date{July 2007}
\begin{document}%

\maketitle

\newtheorem{sub}{\name}[section]
\newtheorem{subn}{\name}
\newtheorem{Thm}{Theorem}[section]
\newtheorem{Lem}[Thm]{Lemma}
\newtheorem{Prop}[Thm]{Proposition}
\newtheorem{Cor}[Thm]{Corollary}
\newtheorem{Rem}[Thm]{Remark}
\newtheorem{Def}[Thm]{Definition}
\newtheorem{Ex}[Thm]{Example}
\renewcommand{\thesubn}{}
\newcommand{\dn}[1]{\def\name{#1}}
\newcommand{\be}{\begin{equation}}
\newcommand{\ee}{\end{equation}}
\newcommand{\bs}{\begin{sub}}
\newcommand{\es}{\end{sub}}
\newcommand{\bsn}{\begin{subn}}
\newcommand{\esn}{\end{subn}}
\newcommand{\bea}{\begin{eqnarray}}
\newcommand{\eea}{\end{eqnarray}}
\newcommand{\BA}[1]{\begin{array}{#1}}
\newcommand{\EA}{\end{array}}
\newcommand{\Real}{\mbox{${\rm I\!R}$}}
\newcommand{\real}{{\rm I\!R}}
\newcommand{\Nat}{\mbox{${\rm I\!N}$}}
\newcommand{\thkl}{\rule[-.5mm]{.3mm}{3mm}}
\newcommand{\Proof}{\mbox{\noindent {\bf Proof} \hspace{2mm}}}
\newcommand{\mbinom}[2]{\left (\!\!{\renewcommand{\arraystretch}{0.5}
 \mbox{$\begin{array}[c]{c}
 #1\\ #2
 \end{array}$}}\!\! \right )}
\newcommand{\brang}[1]{\langle #1 \rangle}
\newcommand{\vstrut}[1]{\rule{0mm}{#1mm}}
\newcommand{\rec}[1]{\frac{1}{#1}}
\newcommand{\set}[1]{\{#1\}}
\newcommand{\dist}[2]{\mbox{\rm dist}\,(#1,#2)}
\newcommand{\opname}[1]{\mbox{\rm #1}\,}
\newcommand{\supp}{\opname{supp}}
\newcommand{\mb}[1]{\;\mbox{ #1 }\;}
\newcommand{\undersym}[2]
 {{\renewcommand{\arraystretch}{0.5}
 \mbox{$\begin{array}[t]{c}
 #1\\ #2
 \end{array}$}}}

\newlength{\wex}  \newlength{\hex}
\newcommand{\understack}[3]{%
 \settowidth{\wex}{\mbox{$#3$}} \settoheight{\hex}{\mbox{$#1$}}
 \hspace{\wex}
 \raisebox{-1.2\hex}{\makebox[-\wex][c]{$#2$}}
 \makebox[\wex][c]{$#1$}
  }%

\newcommand{\smit}[1]{\mbox{\small \it #1}}
\newcommand{\lgit}[1]{\mbox{\large \it #1}}
\newcommand{\scts}[1]{\scriptstyle #1}
\newcommand{\scss}[1]{\scriptscriptstyle #1}
\newcommand{\txts}[1]{\textstyle #1}
\newcommand{\dsps}[1]{\displaystyle #1}

\def\ga{\alpha}     \def\gb{\beta}       \def\gg{\gamma}
\def\gc{\chi}       \def\gd{\delta}      \def\ge{\epsilon}
\def\gth{\theta}                         \def\vge{\varepsilon}
\def\gf{\phi}       \def\vgf{\varphi}    \def\gh{\eta}
\def\gi{\iota}      \def\gk{\kappa}      \def\gl{\lambda}
\def\gm{\mu}        \def\gn{\nu}         \def\gp{\pi}
\def\vgp{\varpi}    \def\gr{\rho}        \def\vgr{\varrho}
\def\gs{\sigma}     \def\vgs{\varsigma}  \def\gt{\tau}
\def\gu{\upsilon}   \def\gv{\vartheta}   \def\gw{\omega}
\def\gx{\xi}        \def\gy{\psi}        \def\gz{\zeta}
\def\Gg{\Gamma}     \def\Gd{\Delta}      \def\Gf{\Phi}
\def\Gth{\Theta}
\def\Gl{\Lambda}    \def\Gs{\Sigma}      \def\Gp{\Pi}
\def\Gw{\Omega}     \def\Gx{\Xi}         \def\Gy{\Psi}

\centerline{\bf Abstract}
\medskip

\small
For a countable-state Markov decision process we introduce an embedding which produces
a finite-state Markov decision process. The finite-state embedded process has the same optimal
cost, and moreover, it has the same dynamics as the original process when restricting to
the approximating set. The embedded process can be used as an approximation which,
being finite, is more convenient for computation and implementation.
\bigskip

\bigskip

\bigskip

\bigskip

\bigskip

\noindent
---------------------------
\bigskip

\noindent

\normalsize

\bigskip

\noindent
{\bf AMS} subject classification: 60G70, 60J10, 62L20.
\newpage
\mysection{Introduction} \label{intro}
In this paper we develop a tool that is useful in studying
countable state Markov Decision Processes (MDPs)~\cite{P}.
A Markov Decision Process is a controlled dynamical system with probabilistic
transitions, that are influenced by the control actions (for precise
definitions see~\S\S~\ref{ss:form}). We consider
discrete-time MDPs with a discrete {\em state space}
$X $
which is either {\em finite} or {\em countably infinite},
to which we will refer in the sequel as {\em countable}.
The cost under consideration is the long-time average cost.

Countable MDPs are obviously more difficult to study, analytically and
numerically, than finite state MDPs.
Several approaches were developed to deal with this issue. The first approach
is to reduce the state space by clustering together ``equivalent" states: see
e.g.~\cite{GDG} and references therein. This approach provides a smaller state
space and exact relations, but requires a very special structure of the MDP in
order for the derived model to have a finite state space. Namely, equivalent
states must have the same transition probability into and out of the state,
under any action, and the same immediate cost.
This special structure seldom exists in applications.

A second approach
approximates countable MDPs by finite
state MDPs using a truncation of the state space.
Existing results show that as the size of the approximating MDP
increases, its cost function and, under some conditions its
optimal policies approach those of
the original, countable MDP. See, e.g.~\cite{C,Ap,Ab} and references therein.
This approach is applicable in greater generality, but typically provides
approximations without an error estimate---thus the results are ``asymptotic"
in nature.

We propose a different approach, with the advantage that
the optimal cost of the
approximating, finite MDP agrees with that of the countable MDP.
Moreover, restricted to the approximating set the optimal policies agree as
well. Thus the term ``exact approximations." 
The main idea is finite embedding. In \cite{F} embedding techniques were used
to obtain optimal policies within various classes.
We apply the embedding approach for approximations with
general, compact action spaces.
Section \ref{sec:appl} develops some applications, where in some cases
exact closed-form expressions can be obtained.

We conclude this section with a precise statement of the problem.
In \S \ref{sec:embed} we introduce the main idea---the finite embedding,
and prove its existence. In \S \ref{sec:cost} we show that the embedding
possesses the desired properties.
We discuss some extensions in \S \ref{sec:ext}.

\subsection{Problem formulation}\label{ss:form}
Consider a process with state space $X \subset \{0,1,2,3,...\}$.
When the system occupies state $i\in X$, then
the controller can influence its behavior by choosing
{\em an action} $a$ from the compact {\em action set} $A_i$, $i\in
X$, which is a subset of the {\em action space} $A$.
Choosing an action $a\in A_i$ has a twofold effect:
%
\begin{enumerate}
\item[(i)]
A {\em running cost} $c(i,a)$ is incurred,
\item[(ii)]
The system transits from state $i$ to $j$
according to the {\em transition probability} $P(j|i,a)$.
\end{enumerate}

\noindent
Thus an MDP is defined in terms of a quadruplet
$${\cal M} = (X,\{A_i\},c(i,a),P(j|i,a)).$$
The state and action at time $k\geq 0$ are denoted $x_k$ and $a_k$
respectively, so that the system's behavior on the infinite time
interval is described in terms of the stochastic process
$\{ (x_k,a_k)\}_{k=0}^{\infty }$.
\smallskip

\noindent
{\bf Admissible policies.} A policy $\pi $ is a law which is used
to choose the actions $a_k\in A_{x_k}$. It is admissible if its choice at
time $k$ depends only on the history
$
(x_0,a_0,...,x_{k-1},a_{k-1},x_k)
$
of the system up to time $k$.
A policy can be either deterministic or randomized,
so that the choice of $a_k$ may be made
according to a probability measure on $A_{x_k}$.
A (possibly randomized) policy which depends only on $x_k$ is called a
``stationary Markov policy''. Such policies generate a state process
$\{x_k\}_{k=0}^{\infty }$ which is a
Markov chain with stationary transition probabilities.
\smallskip

\noindent
{\bf The cost criterion.} An admissible  policy $\pi $ generates the stochastic
processes $\{x_k\}_{k=0}^{\infty }$ and $\{a_k\}_{k=0}^{\infty }$,
and the expected cost flow
\be \label{cstfl}
C_N(i,\pi )=E_i^{\pi }\sum_ {k=0}^{N-1}c(x_k,a_k),\,N\geq 1.
\ee
The expectation
$E_i^{\pi }$ in (\ref{cstfl}) is with respect to the
probability measure $P_i^{\pi }$ induced by $\pi $ on the
set of sequences $\{(x_k,a_k)\}_{k=0}^{\infty }$ with $ x_0 = i $. We address the
optimal control problem of minimizing the functional
$$\pi \mapsto J_i(\pi )=\liminf _{N\to \infty }\frac {1}{N}
C_N(i,\pi ),\;x_0=i,
$$
over all admissible policies,
and call $J_i(\pi )$ {\em the expected long-run average cost}.
The notation
$g^{\star }({\cal M}) $ for the optimal cost
makes explicit the
model ${\cal M}$ under consideration.
An {\em optimal policy} realizes the minimal
long-run average cost, and an $\epsilon$-optimal policy
realizes it up to $\epsilon $.
\smallskip

We need to exclude one case in which it is not possible to approximate a
countable MDP by a finite one.
\begin{Def} \label{drifting}
Let $\sigma $ be a stationary Markov policy of an MDP ${\cal M}$, generating the
state-action process $\{(x_k,a_k)\}$. We say that {\em $\sigma $ is a
drifting policy} if for every finite set $F\subset X$ and every initial
condition $ i \in F $,
\be \label{drftpolcy}
P_i^{\sigma }(x_k\in F)\to 0 \mbox { as } k\to \infty.
\ee
We say that an MDP ${\cal M}$ is {\em a drifting MDP} if there exists a constant
$\delta >0$ such that any stationary Markov policy of ${\cal M}$, say $\sigma $,
that satisfies
\begin{equation}\label{e:drift-cost}
J(\sigma )<g^{\star }({\cal M})+\delta
\end{equation}
is drifting. An MDP is called {\em non-drifting} if it is not
a drifting MDP.
\end{Def}
Obviously it is not possible to approximate a drifting MDP by a finite MDP
while preserving both cost structure, optimality and dynamics. Moreover,
if the MDP is drifting 
then by Definition~\ref{drifting}
for a certain $\delta >0$, for every $0<\nu <\delta $, every
stationary Markov policy with average cost smaller than $g^{\star }+\nu $,
eventually leaves every finite set $F$
with probability 1. But this yields the existence of a Markov policy 
(not necessarily stationary)
with average cost $g^{\star }$ (that is, an optimal policy) which 
eventually leaves every finite set
with probability 1. Assuming that this does not happen 
is the non drifting condition.
In this paper we will therefore consider only non-drifting MDPs.\\
Borkar's coercive condition~\cite{B} requires that the immediate cost is
higher than $g^{\star }$ for ``far states"---as in~(\ref{e:coerc})---and
thus ensures a non-drifting condition.
\mysection{The embedding}\label{sec:embed}
We wish to associate with the given MDP ${\cal M}$ a finite state MDP ${\cal M}_0$,
$$
{\cal M}_0=(X_0,A_0,Q_0(j|i,a),c_0(i,a))
$$
with cost flow $C^0_N (i,\pi )$ defined as in (\ref{cstfl}),
in such a manner that the two MDPs will share common optimality properties
and, in some sense, will share transition probabilities.
To this end we introduce the notion of embedding, whose
exact definition is presented below.
In \S \ref{sec:cost} we prove that this definition implies the desired
exact approximation property.

We will next define the embedding notion, which will be followed by a construction
of an embedding ${\cal M}_0 $.
Denote by $\{(x_k,a_k)\}_{k=0}^{\infty }$ the generic state-action process of
${\cal M} $ and by $\{(\xi _k,\alpha _k)\}_{k=0}^{\infty } $ the generic
state-action process of ${\cal M}_0 $. Given a finite subset $Z\subset X$ define
\begin{equation}\label{e:EtaTau}
\eta = \inf \{ k \geq 0 : x_k \not\in Z\} \qquad
\tau = \inf \{ k > \eta : x_{ k} \in Z\} \ .
\end{equation}
Thus if the initial state is in $Z$ then $\eta $ is the first exit time and
$\tau $ is the first return time, while if the initial state is not in $Z$
then $\eta = 0 $ and again $\tau $ is the first return time.
For a stopping time $\nu > 0 $ define
$$
C_\nu (i,\pi )=E_i^{\pi }\sum_ {k=0}^{\nu -1}c(x_k,a_k)
$$
and note that if $\nu = N$, a deterministic integer, this definition agrees
with (\ref{cstfl}). (If, however, $\nu $ is a random variable then
$ C_\nu (i,\pi ) \not= C_N (i,\pi )\mid_{N= \nu}$,
the latter being a random variable.)
Given a finite subset $Z_0 \subset X_0$,
define $ \eta_0 ,\ \tau_0 $ and $ C_\nu^0 $ analogously.
\begin{Def}\label{embed}
We say that ${\cal M}_0$ is embedded in ${\cal M}$ if there exist subsets
$Z_0\subset X_0$ and $Z\subset X$, and a one-to-one mapping $e:Z_0\mapsto Z$
from $Z_0$ onto $Z$ such that, for any stationary Markov policy $\sigma $
of ${\cal M}$ under which $Z$ contains at least one recurrent state,
the following holds.
\smallskip

\noindent
There exists a stationary Markov policy $\sigma _0$ of ${\cal M}_0$,
such that if the processes start with initial states $ x_0 \in Z $ and
$\xi _0 = e^{-1} (x_0) \in Z_0 $ respectively, then
$x_{\tau }$ and $\xi _{\tau _0}$ are identically distributed
$($under the probability measures $P_{x_0}^{\sigma }$ and
$P_{\xi_0}^{\sigma_0 }$ respectively$)$ and, if the term on the right is finite,
$$C_{\tau }(x_0 ,\sigma )=C^0_{\tau _0}(\xi_0 ,\sigma _0).$$
\end{Def}
Embedding
means that if we restrict attention only to the states $Z$ in $X$ and to the
corresponding states $Z_0$ in $X_0$, then the performance of any stationary
Markov policy $\sigma $ on
${\cal M}$ can be imitated by the performance of some stationary Markov policy $\sigma _0$ on
${\cal M}_0$. This imitation can be achieved when considering finite subsets $F$ of the state
space $X$ on which $\sigma $ generates a nontrivial dynamics, namely the states in $F$ are
not all transient under $\sigma $.

\begin{Rem}
The idea of embedding is a natural generalization of Kac's Theorem and the
``chain on $Z$" idea~\cite[Thm.~10.2.3, \S 10.4.2 and 10.6]{MT}. Kac's Theorem
gives the value of the stationary probability of a Markov chain in terms of
the ``cycle times." Moreover, the stationary distribution of a chain can be
obtained from that of the chain restricted to a subset $Z$, again in terms of
excursion times outside $Z$. In our case
we need to account also for the time and costs accrued during the excursion
outside of $Z$.
\end{Rem}

We now establish that for certain finite sets $Z\subset X$ there exists
a finite state ${\cal M}_0=(X_0,A_0,P_0,c_0)$ and
an embedding $e(\cdot )$ of ${\cal M}_0$ in ${\cal M}$.
This embedding will be useful and significant for stationary
Markov policies $\sigma $ which induce a nontrivial dynamics on $Z$.
The embedding is such that if $e(\cdot )$ is defined on $Z_0$, then
$Z=e(Z_0)$ and
$$\#(X_0)=2\#(Z_0).$$

In \S \ref{sec:cost} we will employ the embedding result to establish
existence and characterize optimal policies for certain countable state MDPs.
As we shall see there, we need the embedding to be such that, in addition to
the cost flow, the expected return times agree, that is,
$E \tau = E \tau_0 $.

\subsection{The existence of embedding}
\begin{Thm}\label{existence}
Let ${\cal M}$ be a countable state MDP and let $Z\subset X$ be a finite set.
Then there exists a finite
state MDP ${\cal M}_0$ with state space $X_0$ and an embedding
$e:Z_0\mapsto Z$ of ${\cal M}_0$ in ${\cal M}$ such that
$$
\#(X_0)=2\#(Z_0).
$$
\end{Thm}
The result has non-trivial content provided that $Z\subset X$
contains a recurrent state of some stationary Markov policy of ${\cal M}$.
\smallskip

\noindent
{\em Proof}: We will define  an MDP ${\cal M}_0$, and will then
establish that it has the properties asserted in the theorem. Denote
$$Z=\{z_1,...,z_n\},$$
let $Z_0$ be a finite set
\be \label{Z0def}
Z_0=\{s_1,...,s_n\}
\ee
and let $e:Z_0\mapsto Z$ be defined by
\be \label{embde}
e(s_i)=z_i, \ i=1,2,...,n.
\ee
With each state $s_i$ in $Z_0$ we associate a state $\omega _i$, and
we then define
\be \label{X0def}
X_0=\{s_1,...,s_n\}\cup \{\omega _1,...,\omega _n\}.
\ee
We have to specify ${\cal M}_0$ as a quadruplet
$$
(X_0,(A_0)_s,c_0(s,a),P_0(s'|s,a))
$$
and define explicitly 
$A_0$, $c_0$ and $P_0$. We define
$A_0(s_i)$ $(=(A_0)_{s_i})$ by
\be \label{A0si}
A_0(s_i)=A(e(s_i)) \mbox { for every } 1\leq i\leq n.
\ee
The definition of the action sets $A_0(\omega _i)$ will be given below.

We next define the transition probabilities
$P_0(s|s_i,a)$, $s\in X_0$, $a\in A_0(s_i)$. First, for
$s_i\in Z_0$
\be \label{P_0sis}
P_0(s|s_i,a)=\left \{
\begin{array}{ll}
P(e(s_j)|e(s_i),a) & \mbox { if } s=s_j\in Z_0\\
1-\sum _{j=1}^nP_0 (s_j|s_i,a) & \mbox { if } s=\omega _i\\
0 & \mbox { if } s = w_j,\ j \not= i. 
\end{array}
\right.
\ee
The corresponding cost is defined by
\be \label{c0def}
c_0(s_i,a)=c(e(s_i) ,a)
\ee
for
$a\in A_0(s_i)=A(e(s_i))$.

We now fix $1\leq i\leq n$ and define the action sets and the transition probabilities
for the state $\omega _i\in X_0\setminus Z_0$. Let $\sigma $ be a
stationary Markov policy of ${\cal M}$ such that $Z$ contains at least
one recurrent state.
Recall the definitions (\ref{e:EtaTau}) and let
\be \label{exitprob}
q_j(\sigma )=P^{\sigma }\{x_{\tau }=z_j \mid x_{\eta - 1} = z_i \} .
\ee
This is the probability that the process $\{x_k\}$ will first enter $Z$
through state $z_j$, conditioned on having left $Z$ from $z_i$ while
employing the action $a\in A(z_i)$ specified by $\sigma $.
{\em The action for state $\omega _i$} induced by $\sigma $ is the collection
of $n+1$ nonnegative numbers
\be \label{omegaction}
\alpha (\sigma )=(\lambda q_1(\sigma ),\ldots ,\lambda q_n(\sigma ), c(\sigma))
\ee
where the constant $\lambda $ satisfies $0<\lambda \leq 1$. These two parameters
$\lambda $ and $c(\sigma )$ which define the action $\alpha (\sigma )$ will be specified below.
The quantities $q_1(\sigma ),...,q_n(\sigma )$ are the probabilities associated
in (\ref{exitprob}) with
the fixed state $z_i$ and the stationary Markov policy $\sigma $. We define the
action set of $\omega _i$ to be 
\be \label{omegaactionset}
A_0(\omega _i)=\bigcup _{\sigma }\alpha (\sigma ),
\ee
where the union is over all the stationary Markov policies
under which $Z$ contains a recurrent state.
Of course, two different policies $\sigma _1$ and $\sigma _2$ may give rise to
the same action, that is  $\alpha (\sigma _1)=\alpha (\sigma _2)$.

We now define the transition probabilities from $\omega _i$. If
$\alpha \in A_0(\omega )$,
then
\begin{equation}\label{e:DefAlpha}
\alpha =
(\alpha_1,\ldots ,\alpha_n ,\alpha_{n+1})=(\lambda q_1,\ldots ,\lambda q_n ,c)
\end{equation}
for some constant $0<\lambda \leq 1$, and it follows that
$
\sum _{j=1}^n\alpha _j=\lambda .
$
We then define $P_0(s|\omega _i,\alpha )$ by
\be \label{P_0omegais}
P_0(s|\omega _i,\alpha )=\left \{
\begin{array}{ll}
\alpha _j & \mbox { if } s=s_j\in Z_0\\
1-\lambda & \mbox { if } s=\omega _i\\
0 & \mbox { if } s\not \in Z_0\cup \{\omega _i\}.
\end{array}
\right.
\ee
Thus for every choice of $0<\lambda \leq 1$, the conditional probability
to enter $Z_0$ through $s_j$, given that the process did enter $Z_0$, is $q_j$,
independent of $\lambda $. However, the value of $\lambda $ determines the expected
time that would elapse until entrance, and we choose $\lambda $ in such manner that
this expected time turns out to be equal to
the corresponding time for ${\cal M}$. Namely, if $\sigma_0 $ is the policy that
uses $ \alpha (\sigma ) $ in state $\omega_i $, then $\lambda $ is chosen such that
\be \label{equaltime}
1 + E^{\sigma_0}_{\omega_i} \tau_0 = E^\sigma_{z_i}  [ \tau \mid \eta = 1 ] \
.
\ee
\begin{Rem}
If $z_i $ is not accessible from the state in $Z$ which is recurrent
under $\sigma $ then we do not need to specify actions for $\omega_i $, while
if it is accessible then necessarily it is recurrent, so that the expectation
on the right-hand side of (\ref{equaltime}) is finite. An appropriate
value of $\lambda $ can clearly be
chosen since for $\lambda = 1 $ the left-hand side equals $1$, while as
$\lambda \to 0 $ this expression diverges.

We note that if it is possible to leave $Z$ from $z_i $ in one step then
$$
E^{\sigma }_{z_i} [ \tau \mid \eta = 1 ]
= 1 + \left [ \sum _{x_j\not \in Z}P(x_j|z_i,\sigma (z_i)) \right ]^{-1}
            \sum _{x_j\not \in Z}P(x_j|z_i,\sigma (z_i))E^{\sigma }_{x_j} \tau.
$$
If, however, this is not possible, then we can define the action in
$A_0(\omega_i )$ that corresponds to $\sigma $ in an arbitrary manner. Finally note
that the normalizing constant $($in square bracket$)$ above is just
$P_0 (\omega_i | s_i , \alpha (\sigma ) ) $.
\end{Rem}
We next consider the cost associated with the action $\alpha $, namely $c_0(\omega _i,\alpha )$,
which is chosen to be such that
\be \label{equalcost}
C_{\tau }(z_i , \sigma ))=C^0_{\tau _0}(e^{-1}(z_i) , \sigma _0))
\ee
holds for all $ z_i \in Z $. We distinguish between two cases:
If under $\sigma $ the process $\{ x_k \} $ starting at $ z_i $ does not leave
$Z$, then equality holds in (\ref{equalcost}) by definition. If the process does
leave $Z$, than
$$
C_\tau ( z_i , \sigma )
=  E_{z_i}^\sigma \sum_{k=0}^{\tau - 1} c (x_k , a_k )
$$
can be expressed in the form
$$
c (z_i , \sigma (z_i ) )
    + \sum _{x_j\not \in Z}P(x_j|z_i,\sigma (z_i)) C_\tau (x_j , \sigma )
    + \sum _{z_j \in Z}P(z_j|z_i,\sigma (z_i)) C_\tau (z_j , \sigma ).
$$
We recall that $c_0(\omega _i,\alpha )$ is actually the ($n+1$)te c mponent of
$\alpha (\sigma )$ in (\ref{omegaction}), denoted $c(\sigma )$. It follows that setting
\begin{equation}\label{e:defOmegaCost}
c_0(\omega _i,\alpha ) =
         [ P ( \omega_i | e^{-1} (z_i) , \sigma_0 ) ]^{-1}
         \sum _{x_j\not \in Z}P(x_j|z_i,\sigma (z_i)) C_\tau (x_j , \sigma )
          [ E_{\omega_i }^{\sigma_0 } \tau_0 ]^{-1}
\end{equation}
ensures the desired equality (\ref{equalcost}).
Thus to each $\sigma $ there corresponds an action $\alpha = \alpha (\sigma )$,
and a cost
$$ c_0 (\omega _i,\alpha ) = \alpha_{n+1} = c(\sigma) $$
which in view of the
explicit expression (\ref{e:defOmegaCost}), indeed depends on the state and action alone.

The definition of ${\cal M}_0$ is thus complete, and it follows from
the definition that ${\cal M}_0$
indeed has the properties asserted in the theorem. $\hfill \Box$
\begin{Rem}
Note that the calculation of $C_\tau $ is very similar to (and is as
complicated as) that of the relative value function in the optimality
equation. More precisely, since we are dealing with a single policy, this is
related to the solution of the Poisson equation: see e.g.~\cite[Thm.~9.5]{MS}.
\end{Rem}
\mysection{Existence of optimal policies}\label{sec:cost}
In order to use the embedding result it is needed that some optimal policy will have a
recurrent state within some finite set.
The following is a simple condition under which
there exists a finite subset $Z$ as required in the theorem presented
in the previous section.
Suppose that we have an estimate
\be \label{gestimate}
g^{\star }({\cal M})<\gamma,
\ee
for some $\gamma $, and moreover, for some ordering of the states $\{ x_j \} $
the following holds:
\be \label{liminf}
\liminf_{j\to \infty }\left\{\min\{ c(x_j,a):a\in A(x_j)\}\right\}>\gamma.
\ee
It clearly follows from  (\ref{gestimate}) and (\ref{liminf}) that if $J(\sigma )< \gamma $,
then some finite set $Z$ contains a recurrent state of $\sigma $.
Such a set is, e.g.,
\begin{equation}
Z = \{ x : \min\{ c(x_j,a):a\in A(x_j)\} < \gamma \} .
\end{equation}
An estimate as in (\ref{gestimate}) does not require a computation of the
optimal policy, but can be
obtained by restricting to a special type of policies.
\smallskip

Fix some state $z$ and denote $s = e^{-1} ( z) $. Define
$\nu = \inf \{ k > 0 : x_k = z \} $, and let $\nu_0 $ be defined analogously.
A sum such as the cycle cost $C_\nu (z ,\sigma ) $ below is well-defined if it
is finite when $c$ is replaced by its absolute value $|c| $.

\begin{Thm}\label{equlcost}
Fix a stationary Markov policy $\sigma $ such that $z\in Z$ is recurrent under
$\sigma $, and let $ {\cal M}_0 $ be an embedding as above.
Let $\sigma _0$ be the stationary Markov policy of ${\cal M}_0$
associated with $\sigma $. If
$$ C_\nu (z ,\sigma )=E_z^\sigma \sum_{k=0}^{\nu - 1} c (x_k , a_k )$$
is well defined then
$$
\lim_{N\to\infty} \frac 1 N C_N (z , \sigma )
= \lim_{N\to\infty} \frac 1 N C^0_N (e^{-1} (z) , \sigma _0) .
$$
\end{Thm}
\begin{Rem}\label{r:inf}
Conditions under which the average cost does not depend on the initial state
are standard, and therefore we shall not elaborate on this point.

The proof applies without change when the condition that $ C_\nu (z ,\sigma )$
is well defined is replaced by the condition that the immediate costs $c (x,a)$
are all nonnegative.
\end{Rem}
{\em Proof:}
We denote $s=e^{-1}(z)$, and by construction $s$ is recurrent under
$\sigma_0 $.
We note that $s_j $ is accessible from $s$
under $\sigma_0 $ if and only if $z_j= e(s_j) $ is accessible from $z$ under
$\sigma $. Therefore we may assume that all states in $Z$ (resp.\ $Z_0 $)
communicate, and ignore transient states. It is also convenient to ignore
(or remove from $Z$) states from which $z$ can be reached only by leaving $Z$.

Let the random times $\nu $ and $\nu _0$ be as in the sentence that precedes
Theorem \ref{equlcost}. Under the recurrence assumption the limit of the
average cost flow exists. It
follows from Theorem 17.2.1 of \cite{MT} that under $\sigma $ the
process pair $ \{ x_k , a_k \} $ possesses an invariant probability measure
which we denote by $\pi $. Moreover, almost surely  under $P_z^\sigma $ we have
$$
\lim_{N\to\infty} \frac 1 N \sum_{k=0}^{N-1} c( x_k , a_k )
=   \frac {E_z^\sigma \sum_{k=0}^{\nu - 1} c (x_k , a_k )}
          {E_z^\sigma \nu }= E_\pi c (x,a),
$$
and analogously for $ {\cal M}_0 $:
$$
\lim_{N\to\infty} \frac 1 N \sum_{k=0}^{N-1} c_0( \xi_k , \alpha_k )
=   \frac {E_s^{\sigma _0}\sum_{k=0}^{\nu _0- 1} c_0( \xi_k , \alpha_k )}
          {E_s^{\sigma _0}\nu _0}= E_{\pi _0}c_0(\xi,\alpha).
$$
In view of our assumptions concerning recurrence and existence of cycle costs,
this implies that
$$
\lim_{N\to\infty} \frac 1 N C_N (z , \sigma )
=  \frac {E_z^\sigma \sum_{k=0}^{\nu - 1} c (x_k , a_k )}
          {E_z^\sigma \nu }
$$
and similarly for $C^0 $.
It therefore suffices to establish that
\be \label{equalratio}
\frac {E_z^\sigma \sum_{k=0}^{\nu - 1} c (x_k , a_k )}
          {E_z^\sigma \nu }=
\frac {E_s^{\sigma _0}\sum_{k=0}^{\nu _0- 1} c_0( \xi_k , \alpha_k )}
          {E_s^{\sigma _0}\nu _0}.
\ee
We note that the numerators in (\ref{equalratio}) are the cycle costs corresponding
to $z$ and $s$, assumed to be well defined.
We first deal with the numerator in the left-hand side of (\ref{equalratio}).
Let $I_{A} $ denote the indicator of the set
$A$, that is
$$
I_{A}(x )=\left \{\begin{array}{lll} 1 & \mbox { if } & x\in A\\
0 & \mbox { if } & x\not\in A.
\end{array}
\right.
$$
Recalling the definitions of $ \eta $ and $\tau $ we have\\
$ \displaystyle{E_z^\sigma \sum_{k=0}^{\nu - 1} c (x_k , a_k )}$
$$
{}\quad = E_z^\sigma \sum_{k=0}^{\min (\nu , \eta) - 1} c (x_k , a_k )
+ E_z^\sigma I_{\{\eta <\nu \}}\left (\sum_{k=\eta }^{\tau - 1} c (x_k , a_k )+
\sum_{k=\tau }^{\nu - 1} c (x_k , a_k )\right ).
$$
By the construction of the embedding,
$$
E_z^\sigma \sum_{k=0}^{\min (\nu , \eta) - 1} c (x_k , a_k )
= E_s^{\sigma_0 } \sum_{k=0}^{\min (\nu_0 , \eta_0) - 1} c (\xi_k , \alpha_k )
$$
since, while $ x_k $ is in $Z$, both transition probabilities and immediate
costs agree.  Also
$$
E_z^\sigma I_{\{\eta <\nu \}}\sum_{k= \eta }^{\tau - 1} c (x_k , a_k )
= E_s^{\sigma_0} I_{\{\eta _0<\nu _0\}}\sum_{k= \eta_0 }^{\tau_0 - 1} c (\xi_k , \alpha_k )
$$
by the definition of the costs $ c_0 (\omega_i , \alpha ) $. Finally, 
using the Markov property,
\be \label{taunucost}
E_z^\sigma I_{\{\eta <\nu \}}\sum_{k=\tau }^{\nu - 1} c (x_k , a_k )
= \sum_{z_j \in Z} P_z^\sigma ( \eta < \nu ,\
                                x_\tau = z_j ) C_\nu (z_j , \sigma ) \ .
\ee
Now write
$$
P_z^\sigma ( \eta < \nu ,\ x_\tau = z_j )
= \sum_{z_i \in Z}
P_z^\sigma ( \eta < \nu ,\ x_\tau = z_j \mid x_{\eta - 1} = z_i )
P_z^\sigma ( x_{\eta - 1} = z_i ) \ .
$$
Recalling that $\nu $ is the return time to state $z$, we
express the first probability on the right-hand side as
$$
P_z^\sigma ( \eta < \nu ,\ x_\tau = z_j \mid x_{\eta - 1} = z_i )
= P_z^\sigma ( x_t \not= z,\ 1 \leq t < \eta ,\
               x_\tau = z_j \mid x_{\eta - 1} = z_i ).
$$
We now observe that the right hand side describes the conditional probability of
two events: one before the conditioning, one after. Since this is a Markov
process we have conditional independence and so
$$
P_z^\sigma ( \eta < \nu ,\ x_\tau = z_j \mid x_{\eta - 1} = z_i )
= P_z^\sigma ( \eta < \nu \mid x_{\eta - 1} = z_i )
\cdot  P_z^\sigma ( x_\tau = z_j \mid x_{\eta - 1} = z_i ) .
$$
It follows from (\ref{taunucost}) and the above computation that\\
$ \displaystyle{E_z^\sigma I_{\{\eta <\nu \}}
                     \sum_{k=\tau }^{\nu - 1} c (x_k , a_k )} $
\be\label{e:messy}
= \sum_{z_j,z_i \in Z}
  P_z^\sigma ( \eta < \nu \mid x_{\eta - 1} = z_i )
       \cdot  P_z^\sigma ( x_\tau = z_j \mid x_{\eta - 1} = z_i )
       \cdot  P_z^\sigma ( x_{\eta - 1} = z_i )
       \cdot  C_\nu (z_j , \sigma ).
\ee
Similarly to (\ref{taunucost}) we have the following expression for $\sigma _0$:
\be \label{taunucost0}
E_s^\sigma I_{\{\eta _0<\nu _0\}}\sum_{k=\tau _0}^{\nu _0- 1} c (\xi _k , \alpha _k )
= \sum_{s_j \in Z_0} P_s^{\sigma _0}( \xi_{\tau _0}= s_j ) C^0_{\nu _0}(s_j , \sigma _0).
\ee
We now repeat the discussion that appears in the text between equations
(\ref{taunucost}) and (\ref{e:messy}) for the embedded process.
Since
all the probabilities and conditional probabilities  in~(\ref{e:messy}) agree
with the corresponding quantities of the embedded process, in view
of (\ref{taunucost}) and (\ref{taunucost0}) it remains to establish that
\be \label{equalCC0}
 C_\nu (z_j , \sigma ) = C^0_{\nu_0 }(s_j , \sigma_0 ).
\ee
  We proceed as before to consider two cases.
We compute $C_{\nu }(z_j,\sigma )$ as follows:\\
$
\displaystyle{C_\nu (z_j , \sigma )
= E_{z_j}^\sigma \sum_{k=0}^{\nu - 1} c (x_k , a_k )}
$
\be \label{eq1}
= E_{z_j}^\sigma \sum_{k=0}^{\min (\nu , \eta) - 1} c (x_k , a_k )
+ E_{z_j}^\sigma \left\{I_{ \{ \eta < \nu \}}
   \sum_{k=\eta}^{\nu - 1} c (x_k , a_k )\right\} .
\ee
The first term once again agrees with the embedded chain. Conditioning
on the exit point we can write the second term as
\be \label{eq2}
E_{z_j}^\sigma \left \{I_{\{ \eta < \nu \}}
   \sum_{k=\eta}^{\nu - 1} c (x_k , a_k )\right \}
= \sum_{z_i \in Z} P_{z_j}^\sigma ( \eta < \nu , x_{\eta - 1} = z_i )
                     C_\nu (z_i , \sigma ),
\ee
and similarly
\be \label{eq3}
E_{s_j}^{\sigma _0}\left \{I_{\{ \eta _0< \nu _0\}}
   \sum_{k=\eta _0}^{\nu _0- 1} c_0 (\xi_k , \alpha _k )\right \}
= \sum_{s_i \in Z_0} P_{s_j}^{\sigma _0}( \eta _0< \nu _0, \xi _{\eta _0- 1} = s_i )
                     C^0_{\nu _0}(s_i , \sigma _0).
\ee
By construction
$$
P_{z_j}^\sigma ( \eta < \nu , x_{\eta - 1} = z_i )
= P_{s_j}^{\sigma_0} ( \eta_0 < \nu_0 , \xi_{\eta_0 - 1} = s_i ),
$$
and since $z$ is accessible from $z_j $, $ P_{z_j}^\sigma ( \eta < \nu ) < 1 $.
Iterating equations (\ref{eq1}), (\ref{eq2}) and (\ref{eq3}) we see that the costs for
both models agree, up to a last term that goes to zero geometrically  fast with the
number of iterations. It follows that the numerators of both sides in (\ref{equalratio})
agree, and the proof for the denominators is
similar. Thus (\ref{equalratio}) is established, and the proof of the theorem is complete.
\hfill $\Box$
\begin{Thm}\label{t:main}
Let $\cal M$ be a Markov Decision Process, and suppose that the state $z$ is
recurrent under a stationary Markov
optimal policy $\sigma $, and that the cycle cost is finite.
Then for any embedding such that $ z \in Z $,
the optimal cost of ${\cal M}_0 $ agrees with that of $\cal M$. Moreover, ${\cal M}_0$
has an optimal policy $\sigma _0$ that agrees with $\sigma $ on corresponding states
of $Z$ and $Z_0$.
\end{Thm}
\noindent
The theorem assumes explicitly that there exists a stationary Markov optimal
policy. This holds for most applications: for conditions see for
example~\cite{P} and references therein.
\mysection{Extensions}\label{sec:ext}
First, note that the requirement that the cycle cost is finite holds whenever
the immediate costs are bounded, since we assume recurrence. Moreover, as
noted in Remark~\ref{r:inf},
this requirement is not needed when the immediate costs are all
of the same sign (positive or negative).

The following result is immediate, but nonetheless useful.
\begin{Thm}\label{t:ActElim}
Fix some $i$.
Suppose the stationary policies $\sigma $ and $\sigma '$ are
such that the associated actions $\alpha (\sigma )$ and $\alpha (\sigma ')$
have costs starting at $\omega_i $ that satisfy
$$ \alpha_{n+1}(\sigma ) =c(\sigma) >c (\sigma ')=\alpha _{n+1}(\sigma '),$$
while
$$ \alpha_i (\sigma ) = \alpha _i (\sigma^\prime ) \mbox { for every } i = 1 , \ldots , n .$$
Then the action $\alpha (\sigma ) $ may be eliminated from $A_0(\omega_i )$.
\end{Thm}
This is quite clear from the definition, and in fact this follows from
standard results of action elimination in MDPs \cite{L}.

Next note that, even if the excursion costs $C_\tau $ are difficult to
calculate, any approximation of $C_\tau $ and of the mean excursion times
leads to a non-exact, approximate embedding, in the sense that optimal costs are not
equal anymore. However, it is easy to see that the approximation is continuous
in the sense that as the approximations of $ C_\tau $ and $E {\tau }$ improve,
the costs (including optimal costs) of the embedded model approach those of the
original MDP.

We now outline the extension to constrained MDPs, where a detailed description
of the model may be found in \cite{Ab}.
In addition to the usual four components of an MDP we define a collection of
immediate cost functions $\{ d_k (x,a),\ k= 1 , \ldots , K \} $. Define
$J^k_i (\pi ) $ in the same way that $J_i (\pi ) $ is defined, but with $d_k $
replacing the immediate cost $c$. The constrained optimization problem is
to minimize the functional $J_i (\pi ) $, subject to the constraints
$$
J_i^k (\pi ) \leq V_k
$$
for some prescribed constants $V_k $, $1\leq k\leq K$.

Standard approximations of constrained problems are more difficult to handle
and establish than unconstrained approximations of
optimization problem. The reason for this is that when we
require the approximate model to satisfy the hard constraints
$$
J_i^{0,k} (\pi ) \leq V_k,\,k=1,2,...,K,
$$
then clearly we may lose continuity, in the sense that even if the original
problem is feasible (that is, there exist policies satisfying the
constraints), an approximation of the required type may not be feasible~\cite{Ab}.
However, using our exact approximation, this difficulty does not arise.

The embedding results hold for this model, with the following minor
modification. Recall the definitions~(\ref{e:DefAlpha})
and~(\ref{e:defOmegaCost})
of the action in ${\cal M}_0 $. Define $d^k (\sigma ) $ as
in~(\ref{e:defOmegaCost}) and define $\alpha (\sigma ) $ by
\begin{equation}
\alpha =
(\alpha_1,\ldots ,\alpha_n ,\alpha_{n+1},\alpha_{n+2},\ldots ,\alpha_{n+K+1})
=(\lambda q_1,\ldots ,\lambda q_n ,c , d^1 , \ldots , d^K) .
\end{equation}
Then the same arguments show that for the embedded chain, all costs agree with
those of the original model, so that we may approximate the countable chain by
a finite chain.

Note, however, that this model is much less robust. Whereas small errors in
the calculation of the cycle cost for the optimization problem may, in the
worst case, lead to sub-optimality, in the constrained case such errors may
lead to infeasibility, that is, violation of the constraints. Thus, if
cycle costs can not be computed exactly, this approximation shares the
infeasibility problem with other, more traditional approximation methods.
\noindent
\mysection{Examples}\label{sec:appl}
\begin{Ex}{\rm
Markov Decision Processes serve as a common model in the control of
dams and reservoirs~\cite{LB}. 
It is standard to discretize both space and time, in order to
arrive at a manageable model. So, let us model the inflow into the water 
reservoir by a Markov chain $D_t $ and let $L_t $ denote the water level at the
reservoir. At each epoch (usually month), the decisions are to use $Z_t $
of the water for electricity generation, and evacuate $Y_t $ through spillways.
The water level is then given by
$$
L_{t+1} = L_t + D_t - (Y_t + Z_t ) .
$$
Denote $x_t = (L_t , D_t ) $ and $a_t = ( Y_t , Z_t ) $. The revenue (negative
cost) in the model arises from the sale of electrical power. 
The amount of power produced is a linear function of the amount of water $Z_t $
used for generation, but also depends non-linearly on the water level $L_t $ 
(since higher water level entails higher energy per unit of water). Thus the
cost $c = c (l,z) $ depends on the state and action taken.
Since $D_t $
is a Markov chain, we obtain an MDP
with discrete state process $x_t $ and discrete action process $a_t $ 
(assuming that the discretization of $D_t $
is compatible with that of $Y_t $ and $Z_t $).
Both control variables are positive and
bounded due to practical considerations (limited
capacity of the generators, and of the spill mechanisms),
so that the number of control actions is finite. Denote the maximal allowed
value of $Y_t $ (resp.~$Z_t $) by $\bar y $ (resp.~$\bar z $).

Since water levels which are too high may pose danger, if water level and
inflow rates are too high, the maximal value of the control variables must be
used, that is $ Y_t = \bar y $ and $Z_t = \bar z $.
This may be formulated by defining a finite set $F \subset X $ so that if the
state $ (l,d) $ is outside this set, the allowed action is only $(\bar y, \bar
z )$. It is therefore natural to use our embedding results so that only states
in $F$ need be considered. We note that since the state space is two
dimensional, the number of states that we ignore (outside $F$) may be
significant.

The case where $D_t $ is i.i.d.\ is particularly simple, since in this case
the state process is one-dimensional, and outside $F$ the process behaves
exactly like a random walk. Therefore, results on random walks may be used to
calculate the entrance distribution and mean time. In fact, in the simple case
that the cost satisfies $c (l,\bar z ) = \bar c $, a constant, for water level
outside of $F$, the cycle cost is just a constant multiple of the average
return time. For an explicit calculation in a simple case, see the next
example.

In water reservoir applications, the object for optimization is often a group
of reservoirs: there could be dozens of reservoirs, all connected. In this
case the ``curse of dimensionality" makes it impossible to solve such models,
and additional approximations are required. For a step in this direction see
the multi-dimensional queueing problem below.
}
\end{Ex}

\begin{Ex}{\rm
In the simplest case where $D_t $ are Bernoulli, and where the release $Y_t +
Z_t $ can only take the values $0 $ or $1$ (with some probability which we can
choose),
we arrive at a generic model of operations research.
Consider the problem of controlling a single queue. New jobs arrive according
to an i.i.d.\ sequence of Bernoulli random variables with mean $\lambda $, and
join an infinite queue. The job at the head of the queue is served, and
(independently) the
probability of completion of service (representing the speed of service) is
the control variable $a$. Assume that there exists some $I> 0 $ and
$ \mu > \lambda > \delta $ so that
$$
A (x) = [\delta , \mu ] \ \mbox{for $x < I $} \quad \mbox{and}\
A (x) = \{ \mu \} \ \mbox{for $x \geq I $} .
$$
That is, service rate is controlled for small queue size, but maximal rate
must be used if the queue is large. If we assume further that, for $x \geq I$
the immediate cost is sub-linear, that is, it
satisfies $ c (x,a) = c (x , \mu ) \leq c \cdot x $ for some positive $c$,
then all our assumptions hold: state $0$ is recurrent under any policy and
cycle costs are finite.
The embedded model ${\cal M}_0 $ can be computed as follows.
Set $X_0 = \{ 0 , 1 , \ldots , I - 1 \} $. Since the process moves at most one
step per unit time, the only way to exit this set is through state $I$, so we
can set $ \omega = I $. To calculate $c_0 $, we need
$$
 C_\tau (I , \sigma ) = E_I \sum_{k=0}^{\tau-1} c (x_k , \mu )
$$
which in fact does not depend on $\sigma $, since the only available action is
$\mu $. This expression is identical to a standard queue, without control,
with arrival rates $\lambda $ and service rate $\mu $.
However, for the standard queue, transition probabilities for $ x \not= 0 $
do not depend on the state $x$. So, denote by $E^0$ expectation for the
standard queue, and we can calculate
$$
C_\tau (I , \sigma ) = E^0_0 \sum_{k=0}^{\tau-1} c (x_k - I , \mu ) \ .
$$
But this equals $ E_0^0 \tau \cdot J $ where $J$ is the average cost for a
standard queue. This can be calculated in terms of the stationary distribution
$\pi $ of the standard queue, that is
$$
C_\tau (I , \sigma ) = [ E_0^0 \tau ]^{-1}\sum_{x=0}^\infty \pi_x c (x , \mu )
\ . 
$$
Define $ \rho = \lambda ( 1 - \mu ) \slash \mu ( 1 - \lambda ) $. Then it is
easy to check that for the standard queue, $\pi_x = (1-\rho ) \rho^ x $.
On the other hand, by Kac's Theorem~\cite[Thm.~10.2.3]{MT},
$ E_0^0 \tau = [\pi_0 ]^{-1} = 1 \slash ( 1 - \rho )$
so that we have an explicit expression for $ C_\tau (I , \sigma ) $.
Clearly, this can be calculated analytically for various costs functions
$c$, and in particular for linear costs.

This example can be extended as follows. Suppose we do not assume that the
action space is restricted for $ x \geq I $. Instead assume that
\begin{equation}\label{e:coerc}
c (x,a) < \min \{ c ( y,a ) : a \in A_y ,\ y \geq I \}
\end{equation}
for all $ x < I $ and all $a$. It follows that $ \mu $ is the
optimal action at $ x \geq I $, and we can apply Theorem~\ref{t:ActElim}
so that the previous conclusions apply.

In general, since this is a skip-free,
one-dimensional problem, our results allows an easy
decoupling of the behavior for $ x < I $ from that for $ x \geq I $. The
situation is more complicated if the skip-free assumption is violated,
namely either batch arrivals or batch service or both are allowed.
However, as is clear from the proof of Theorem~\ref{t:main}, we can write an
implicit expression for the cost using the cost flows until the first hitting
time of $ \{ 0 , 1 , \ldots , I- 1 \} $, to obtain an explicit expression for the
embedded model.
}
\end{Ex}

\begin{Ex}
{\rm
Consider now a multi-dimensional queueing problem. Jobs of type $1,\ldots ,K $
arrive according to a $K$-dimensional process $B(t)$
of i.i.d.\ vectors. The $k$th
coordinate represents arrivals of customers of type $k$. Customers join
infinite queues, one queue for each type. A single server chooses at each
point in time which queue to serve, and serves the job at the head of the
line. If job of type $k$ is served, then the service will succeed
with probability $\mu_k $, and then the job will leave the queue.

If we impose the condition that some queue must be served as long as not all queues
are empty, then the empty state will be recurrent under mild conditions. For
example, it is sufficient to assume that the total number of arrivals at any unit time
interval is bounded, and that the condition
$$
\sum_{k=1}^K \frac {E B_k }{\mu_k} < 1 \ .
$$
Let $Q(t)$ be the vector of queue-sizes at time $t$. This is the state of our
MDP, and the control is the choice of queue to serve. This is easily seen to
be an MDP, once immediate costs $c (q,a) $ are chosen.
It is then natural to choose
$$
Z = \left \{ Q : \sum_{k=1}^K Q_k (t) \geq Q_0 \right \}
$$
and approximate the infinite model with a finite one.

Suppose for some $Q_0 $ we have that if
$$
\sum_{k=1}^K Q_k (t) \geq Q_0
$$
then
$$
c (x,a) = \sum_{k=1}^K c_k Q_k
$$
for some positive coefficients $\{ c_k \} $.
Then the system simplifies considerably, and the computation of the hitting
distributions and costs, required for our approximation, become
feasible~\cite{BMM,W}.
}
\end{Ex}

\bigskip
\noindent{\bf Acknowledgements:}

\noindent
We would like to thank the reviewers for their insightful comments.

\noindent
This research was supported in part by the fund for promotion of research at
the Technion and the fund for promotion of sponsored research at
the Technion.\\
Adam Shwartz holds the Julius M.\ and Bernice Naiman Chair in Engineering at
the Technion.

\end{document}